\documentclass[12pt,twoside]{article}

\usepackage{amssymb,amsmath,amsthm}
\usepackage{times}

\newtheorem{theorem}{Theorem}
\newtheorem{lemma}[theorem]{Lemma}
\newtheorem{proposition}[theorem]{Proposition}
\newtheorem{remarka}[theorem]{Remark}

\newtheorem{definitiona}[theorem]{Definition}

\newcommand{\e}{\mathrm{e}}
\newcommand{\dd}{\mathrm{d}}
\newcommand{\ii}{\mathrm{i}}
\newcommand{\R}{\mathbb{R}}

\renewcommand{\Re}{\operatorname{Re}}

\newcommand{\bs}{\backslash}
\newcommand{\half}{\mbox{$\frac{1}{2}$}}

\newcommand{\pdpd}[2]{\frac{\partial #1}{\partial #2}}
\newcommand{\dede}[2]{\frac{{\mathrm{d}} #1}{{\mathrm{d}} #2}}

\newcommand{\supp}{\operatorname{supp}}

\newcommand{\Op}{\operatorname{Op}}

\newcommand{\Hamx}{\Gamma_{x}}
\newcommand{\Hamxi}{\Gamma_{\xi}}

\newcommand{\Bsub}{B_{\rm s}}
\newcommand{\alphatilde}{{\widetilde{\alpha}}}
\newcommand{\What}{{\widehat{W}}}

\begin{document}
\pagestyle{myheadings}
\markboth{Christiaan C. Stolk}{%
A hyperbolic initial value problem with dissipation}

\title{Parametrix for a hyperbolic initial value problem with dissipation 
in some region}
\author{Christiaan C. Stolk\\
        Centre de Math\'ematiques\\
        Ecole Polytechnique\\
        91128 Palaiseau Cedex\\
        stolk@math.polytechnique.fr}
\date{December 5, 2003}
\maketitle

\begin{abstract}
\noindent
We consider the initial value problem for a pseudodifferential
equation with first order hyperbolic part, and an order $\gamma > 0$
dissipative term. Under an assumption, depending on an integer
parameter $L \geq 2$ such that $2 \gamma < L$, we construct for this
initial value problem a parametrix that is a Fourier integral
operator of type $\rho = 1 - \gamma/L$. The assumption implies that
where the principal symbol of the dissipative term is zero, the
terms of order up to $L-1$ in its Taylor series also vanish.
\end{abstract}



\section{Introduction}

In this paper we study pseudodifferential operators of the form
\begin{equation} \label{eq:define_P}
  P = \partial_z - \ii A(z,x,D_x) + B(z,x,D_x) .
\end{equation}
We assume that $A$ and $B$ satisfy\\[.5ex]
(i) $A = A(z,x,D_x)$ is a smooth family of pseudodifferential operators
in $\Op S^1 (\R^n \times \R^n)$, with homogeneous real principal
symbol $a=a(z,x,\xi)$.\\[.5ex]
(ii) $B = B(z,x,D_x)$ is a smooth family of pseudodifferential operator
in $\Op S^\gamma(\R^n \times \R^n)$, $\gamma > 0$ with non-negative
real homogeneous principal symbol $b=b(z,x,\xi)$. To be precise we
assume $b$ homogeneous for $|\xi| \geq 1$, and smooth for
$|\xi| < 1$.\\[.5ex]
(iii) The derivatives up to order $L-1$, $L \geq 2$ of $b$ and of
$\Bsub = B - b$ satisfy
\begin{align}
  & \big| \partial_{(z,x)}^\alpha \partial_\xi^\beta b(z,x,\xi) \big|
  \leq C (1+|\xi|)^{-|\beta|+\frac{|\alpha|+|\beta|}{L}\gamma}
    (1+b(z,x,\xi))^{1-\frac{|\alpha|+|\beta|}{L}} , \quad\!
    |\alpha| + |\beta| < L ,
    \label{eq:assumption_b1}\\
  & \big| \partial_{(z,x)}^\alpha \partial_\xi^\beta \Bsub(z,x,\xi) \big|
  \leq C (1+|\xi|)^{-|\beta|-1+\frac{|\alpha|+|\beta|+2}{L}\gamma}
    (1+b(z,x,\xi))^{1-\frac{|\alpha|+|\beta|+2}{L}} ,
    \label{eq:assumption_b2} \\
  & \makebox[11cm]{}|\alpha| + |\beta| + 2 < L \nonumber.
\end{align}
Without $B$ the operator $P$ would be a standard first order hyperbolic
operator. The operator $B$ introduces damping into the equation.
Because $B$ is of order $\gamma > 0$, this will lead to the
suppression of singularities propagating in a region where $b \neq 0$.
The third assumption means that $B$ increases slowly near points where
$b = 0$. It implies that the derivatives of $b$ up to order $L-1$
vanish on the set of $(z,x,\xi)$ where $b = 0$. It is automatically
satisfied for $L = 2$, $\gamma =1$ by the inequality $b'{}^2 \leq
2 b \| b'' \|_{L^\infty}$.  By (ii) it is also satisfied for $|\alpha| +
|\beta| \geq L$, with a constant $C_{\alpha,\beta}$ depending on
$\alpha,\beta$. We study the initial value problem for $P$ with
initial value at some $z_0 \in \R$, given by
\begin{equation}
  P u(z,\cdot) = 0 , \quad \text{ for } z \in ]z_0,Z[ , \qquad
  u(z_0,\cdot) = u_0 .
    \label{eq:initial_value_problem}
\end{equation}

Operators of this type appear in one-way wave equations, see
\cite{StolkOneWayAcoustic1}. Here $B$ is assumed to vanish on a
large subset of $(\R \times \R^n) \times \R^n$, that has a
non-empty interior. Outside this set the equation should not admit
the propagation of singularities and $B$ is assumed to be
non-zero. As a consequence $B$ vanishes to all orders at the
boundary of the interior of the zero set. Operators of the type
(\ref{eq:define_P}) also appear in drift diffusion
equations. Estimates for exponential decay in hyperbolic first order
systems, with $B = 0$, were given in \cite{RauchTaylor1974}.

The existence and uniqueness of solutions to
(\ref{eq:initial_value_problem}) follows from the results in
section 23.1 of \cite{hormander85a} for $\gamma \leq 1$. In
section~\ref{sec:wellposedness} we extend this to the case where
$B$ satisfies (\ref{eq:assumption_b1}), (\ref{eq:assumption_b2})
with $2\gamma < L$. The solution operator will be denoted
$E : H^s(\R^n) \rightarrow C([z_0, Z];H^s(\R^n))$. By $E(z,z_0)$
we will denote the map $u_0 \rightarrow E u_0 (z,\cdot)$.

The solutions will be related to those of the purely hyperbolic
operator $P_0$ defined by
\begin{equation}
  P_0 = \partial_z - \ii A(z,x,D_x) .
\end{equation}
By $E_0$ and $E_0(z,z_0)$ we denote the solution operator to the
initial value problem for $P_0$. It is well known that $E_0$ is a
Fourier integral operator, and can also be defined for $z \leq
z_0$ if $A$ and $B$ are defined there, see \cite{duistermaat72}
or a text on Fourier integral operators such as
\cite[section 5.1]{duistermaat96} (note that
$A(z,x,D_x)$ is not strictly a pseudodifferential operator in
$(z,x)$, but a $z$-family of pseudodifferential operators in $x$,
but the argument remains valid, see e.g.\ theorem 18.1.35 of
\cite{hormander85a}.) Let $p_0 = p_0(z,x,\zeta,\xi) = \ii \zeta -
\ii a(z,x,\xi)$ denote the principal symbol of $P_0$. For $B=0$
the singularities (elements of the wave front set of the solution)
of the solution propagate on the set $p_0 = 0$ according to the
Hamilton vector field of $p_0/\ii$, which reads
\begin{equation}
  \pdpd{}{z} - \pdpd{a}{\xi} \cdot \pdpd{}{x}
    + \pdpd{a}{x} \cdot \pdpd{}{\xi} .
\end{equation}
The solution curves to this field are called bicharacteristics.
They can be parameterized by $z$. We denote the $(x,\xi)$
components of the solution curve with initial values $(x_0,\xi_0)$
at $z_0$ by
$(\Hamx(z,z_0,x_0,\xi_0),\linebreak[2]\Hamxi(z,z_0,x_0,\xi_0))$.

Let $I = I(z,x,\xi)$ be the integral of $B$ along a
bicharacteristic with initial values $x,\xi$ at $z$ (not at $z_0$)
\begin{equation} \label{eq:define_I}
  I(z,x,\xi) = \int_{z_0}^z
    b(z',\Hamx(z',z,x,\xi),\Hamxi(z',z,x,\xi)) \, \dd z' .
\end{equation}
Let $\Phi_{z,z_0}$ denote the bicharacteristic flow on $\R^n
\times \R^n \bs 0$, i.e.\
\begin{equation}
  \Phi_{z,z_0}(x,\xi) = (\Hamx(z,z_0,x,\xi),\Hamxi(z,z_0,x,\xi)) .
\end{equation}
We also define $\widetilde{I} = \widetilde{I}(z,x,\xi)$ by
\begin{equation}
  \widetilde{I}(z,x,\xi)
  = (I \circ \Phi_{z,z_0})(z,x,\xi)
  = \int_{z_0}^z b(z',\Hamx(z',z_0,x,\xi),\Hamxi(z',z_0,x,\xi)) \, \dd z' .
\end{equation}
The factors $\exp(-I(z,x,\xi))$, $\exp(-\widetilde{I}(z,x,\xi))$
will be used in the construction of solutions to
(\ref{eq:initial_value_problem}).

We will construct a parametrix for (\ref{eq:initial_value_problem}),
in the form of a family of pseudodifferential operators if
$A = 0$, or a Fourier integral operator with real phase function if
$A \neq 0$. The factor $\exp(-I)$ will be part of the
amplitude. A complication is that $\exp(-I)$ is not a standard
symbol in $S^0$. Instead we need the symbol classes
$S^\mu_{\rho,\delta} (\R^n \times \R^n)$ of type $\rho,\delta$. By
definition a function $f$ in $C^\infty(\R^n \times \R^n)$ is in
the space $S^\mu_{\rho,\delta}(\R^n \times \R^n)$ of symbols of
order $\mu$ and type $\rho,\delta$, if there are constants
$C_{\alpha,\beta}$ such that
\begin{equation} \label{eq:symbol_estimate_rho}
  \big| \partial_x^\alpha \partial_\xi^\beta f(x,\xi) \big|
  \leq C_{\alpha,\beta}(1+|\xi|)^{\mu+|\alpha| \delta - |\beta| \rho} .
\end{equation}
The case $\rho = 1,\delta=0$ is the standard case. For $\gamma =
1$ and without additional assumptions on $B$ the problem
(\ref{eq:initial_value_problem}) leads
to operators of type $(\rho,\delta) = (\half,\half)$. For example,
when $A = 0$ and $B$ is independent of $z$ it is known that
$E(z,z_0) = \exp((z-z_0) B(z,x,D_x))$ is in $\Op
S^0_{\frac{1}{2},\frac{1}{2}}$ (see the remarks in \cite[page
515]{treves80b}). This case has also been analyzed using Fourier
integral operators with complex-valued phase function
\cite{Kucherenko1973,melin76}.

Symbols of type $\rho, \delta$ do not always lead to a good
calculus of pseudodifferential and Fourier integral operators. For
pseudodifferential operators many of the usual results hold when
$\rho - \delta > 0$. If on the other hand $\rho - \delta = 0$ then
the commutator of two pseudodifferential operators is no longer of
lower order. In that case a more refined analysis is required,
which we will not discuss here. For Fourier integral operators one
has $\rho = 1 - \delta$, and the requirement becomes $\half < \rho
\leq 1$, see e.g.\ Theorem 2.4.1 in \cite{duistermaat96}.

Here we show that when (\ref{eq:assumption_b1}), (\ref{eq:assumption_b2})
are satisfied with $2 \gamma < L$, then we have symbols in classes
$S^\mu_{1-\frac{\gamma}{L},\frac{\gamma}{L}}$, so that a good
calculus exists. A first indication of such behavior is that, if
$f$ is a symbol $[0,\infty[ \rightarrow \R$, i.e.\ satisfies
$|f^{(k)}(y)| < C_k (1+y)^{\delta - k}$, then $f \circ b$ is a
symbol in
$S^{\max(\delta,0)}_{1-\frac{\gamma}{L},\frac{\gamma}{L}}$. The
proof (almost immediate) is given in
section~\ref{sec:wellposedness}, where we also show that a
pseudodifferential square root $\sqrt{1 + B}$ exists modulo a
regularizing operator. Our first main result is the following
theorem about the solution operator $E$ to the initial value
problem (\ref{eq:initial_value_problem}). When we write that $W$ is a
bounded family of pseudodifferential operators in
$S^{\mu_0}_{\rho,\delta}$ with $\partial_z^j W$ in
$S^{\mu_j}_{\rho,\delta}$, then it will be understood that
$\partial_z^j W$ is also a bounded family in this class.

\begin{theorem} \label{th:main_theorem}
Let $A$ and $B$ satisfy (i), (ii) and (iii), and assume that
$2\gamma < L$. Then there are bounded families of
pseudodifferential operators $W = W(z,x,D_x)$ and $\widetilde{W} =
\widetilde{W}(z,x,D_x)$ in $\Op
S^0_{1-\frac{\gamma}{L},\frac{\gamma}{L}}(\R^n \times \R^n)$,
with $\partial_z^j W$ and $\partial_z^j \widetilde{W}$ in $\Op S^{j \gamma}_%
{1-\frac{\gamma}{L},\frac{\gamma}{L}}(\R^n \times \R^n)$, such
that
\begin{equation} \label{eq:parametrix_left_and_right}
  W(z,x,D_x) E_0(z,z_0)
  = E(z,z_0)
  = E_0(z,z_0) \widetilde{W}(z,x,D_x) .
\end{equation}
The symbol $W(z,x,\xi)$ can be written as an asymptotic sum
\begin{equation} \label{eq:form_of_W_left}
  W = (1 + \sum_{j=1}^\infty K^{(j)}) \exp(-I) + R
\end{equation}
where $R$ is a smooth family of symbols in $S^{-\infty}(\R^n
\times \R^n)$, and the $K^{(j)}$ satisfy
\begin{align}
  & \begin{array}{l}
    K^{(j)} \text{ is a smooth family of symbols in }
    S^{j(2\gamma-1)}(\R^n \times \R^n) , z \in [z_0,Z] , \text{ and}
    \label{eq:Wprop1}
    \end{array} \\
  & \begin{array}{l}
    K^{(j)}\exp(-I) \text{ is a bounded family of symbols in }
    S^{-j(1-\frac{2\gamma}{L})}(\R^n \times \R^n) , z \in [z_0,Z] ,\\
    \text{with } \partial_z^k (K^{(j)}\exp(-I)) \text{ in }
    S^{k\gamma - j(1-\frac{2\gamma}{L})}(\R^n \times \R^n) .
    \end{array}\!\!\!
    \label{eq:Wprop2}
\end{align}
Furthermore, $K^{(j)}(z,x,\xi) = 0$ for $(z,x,\xi)$ such
that $I(z,x,\xi) = 0$. The same is true for
$\widetilde{W}(z,x,\xi)$ with $\widetilde{I}$ instead of $I$.
\end{theorem}

When $I(z,x,\xi) > 0$ for some point $(z,x,\xi)$, then
$\exp(-I(z,x,\lambda \xi))$ and its derivatives decay
exponentially for $\lambda \rightarrow \infty$, and $\exp(-I)$ is
in $S^{-\infty}$. Hence for a bicharacteristic with some finite
part in the region of cotangent space where $B > 0$, the factor
$W$ is in $S^{-\infty}$, and the solution operator becomes
regularizing. The difficult part is therefore the behavior near
the boundary of the region $B>0$.

For each $z$ derivative our estimate of $W(z,x,\xi)$ worsens by a factor
$(1+|\xi|)^\gamma$. It turns out that in fact there is a better
estimate. Let $C_0$ denote the canonical relation of $E$
\begin{multline}
  C_0 = \{ ( z, \Hamx(z,z_0,x_0,\xi_0),
        a(z_0,x_0,\xi_0), \Hamxi(z,z_0,x_0,\xi_0) ; x_0,\xi_0) \, | \, \\
        (x_0,\xi_0) \in \R^n \times \R^n \bs 0, z > z_0 \} .
\end{multline}

\begin{theorem} \label{th:FIO_property}
The symbols $W$ and $\widetilde{W}$ are in $S^0_{1-\frac{\gamma}{L},%
\frac{\gamma}{L}}((]z_0,Z[ \times \R^n) \times \R^n)$, and the
map $E$ is a Fourier integral operator in
$I_{1-\frac{\gamma}{L}}^{-1/4}(]z_0,Z[ \times \R^n, \R^n; C_0)$.
\end{theorem}

Recall that a symbol in $S^\mu_{\rho,\delta}((]z_0,Z[ \times \R^n) \times
\R^n)$ is assumed to satisfy estimates like (\ref{eq:symbol_estimate_rho})
only locally. A function $f=f(z,x,\xi)$ is in $S^\mu_{\rho,\delta}
((]z\-0,Z[ \times \R^n) \times \R^n)$ if for each $\alpha,\beta$ and each
conically compact subset $K$ of $(]z_0,Z[ \times \R^n) \times \R^n)$ there
is a constant $C_{\alpha,\beta,K}$ such that
\begin{equation} \label{eq:symbol_estimate_z_x}
  \big| \partial_{z,x}^\alpha \partial_\xi^\beta f(z,x,\xi) \big|
  \leq C_{\alpha,\beta,K}(1+|\xi|)^{\mu+|\alpha| \delta - |\beta| \rho}
  , \quad (z,x,\xi) \in K .
\end{equation}
There need not be constants such that (\ref{eq:symbol_estimate_z_x}) is valid
globally. Indeed we find bounds for the $z$-derivatives of
$W(z,x,\xi)$ that blow up when $z \rightarrow z_0$.

The organization of the paper is as follows. In the next section
we discuss the well-posedness of the problem
(\ref{eq:initial_value_problem}). We then study the case $A =
0$.
In section~\ref{sec:symbol_estimates} we estimate various
quantities in terms of powers of $(1+|\xi|)$ and $(1+I)$, and we
discuss symbols estimates for functions with a factor $\exp(-I)$
that result from this. We then prove Theorem~\ref{th:main_theorem}
for the case $A = 0$ in section~\ref{sec:parametrix_A=0}. In
section~\ref{sec:complete_result} we use this result and Egorov's
theorem to prove the case $A \neq 0$. We then prove the Fourier
integral operator property in section~\ref{sec:FIO_property}. In
the final section we discuss a choice of $b$ that satisfies the
assumption (\ref{eq:assumption_b1}).

To denote a constant we will use the letter $C$. The value of $C$
can change between equations.

\section{Well-posedness\label{sec:wellposedness}}

In this section we construct the square root $\sqrt{1 + B}$, modulo a
smoothing operator. Using this square root we show that the problem
(\ref{eq:initial_value_problem}) has unique solutions satisfying energy
estimates. The following lemma concerns the square root.

However, we first show that if $f : [0,\infty [ \rightarrow \R$ is a symbol
of order $\delta$, that is for each $k$ there is a $C_k$ such that
$|f^{(k)}(y)| < C_k (1+y)^{\delta - k}$, then $f \circ b$ is a symbol
in $S^{\max(\delta,0)}_{1-\frac{\gamma}{L},\frac{\gamma}{L}}$.
A simple computations shows that
$\partial_{(z,x)}^\alpha \partial_\xi^{\beta} f(b)$ is a sum of
terms of the form
\begin{equation}
  c f^{(k)} \prod_{j=1}^k \partial_{(z,x)}^{\alpha_j}
    \partial_\xi^{\beta_j} b ,
\end{equation}
where $c$ is a constant, $\sum_j \alpha_j = \alpha, \sum_j \beta_j =
\beta$ and $(\alpha_j,\beta_j) \neq 0$. By (\ref{eq:assumption_b1}) this is
less or equal than
\begin{equation}
  C'_{K,\alpha,\beta} (1+|\xi|)^{-\beta + \frac{|\alpha| + |\beta|}{L} \gamma}
    (1+b)^{\delta-\frac{|\alpha| + |\beta|}{L}}
  \leq
  C_{K,\alpha,\beta} (1+|\xi|)^{\max(\delta,0)-\beta
          + \frac{|\alpha|+|\beta|}{L} \gamma} ,
\end{equation}
hence $f \circ b$ is in
$S^{\max(\delta,0)}_{1-\frac{\gamma}{L},\frac{\gamma}{L}}$.

\begin{lemma} \label{lem:root_B}
Assume that $B$ is selfadjoint and satisfies (ii) and (iii) and that
$2\gamma < L$. Then there is a smooth family of pseudodifferential
operators $Q(z,x,D_x) \in
\Op S^{\gamma/2}_{1-\frac{\gamma}{L},\frac{\gamma}{L}}
(\R^n \times \R^n)$ with $\partial_z^j Q(z,x,D_x) \in
\Op S^{\gamma/2+\frac{j\gamma}{L}}_{1-\frac{\gamma}{L},\frac{\gamma}{L}}
(\R^n \times \R^n)$, such that $Q$ is selfadjoint and
$Q^2 = 1 + B + R$, with $R \in S^{-\infty}$.
\end{lemma}

\begin{proof}[Proof of Lemma~\ref{lem:root_B}]
The proof is a variant of the standard construction as an asymptotic
sum (see e.g.\ Lemma II.6.2 in \cite{taylor81})
\begin{equation}
  Q = Q^{(0)} + Q^{(1)} + \ldots .
\end{equation}
In this construction the first term satisfies
$Q^{(0)} = (1+b)^{1/2} + \text{l.o.t.}$ A remainder is defined
\begin{equation}
  R^{(k)} = \left( \sum_{j=0}^{k-1} Q^{(j)} \right)^2 - (B + 1) ,
\end{equation}
and this is used to define the symbol of the next term, such that
\begin{equation}
  Q^{(k)} = - \half (1+b)^{-\frac{1}{2}} R^{(k)} + \text{l.o.t.}
\end{equation}
We must show that the $Q^{(k)}$ defined in this way are symbols of
decreasing order.

As an induction hypothesis we assume that $R^{(k)}$ is selfadjoint and
a sum of terms of the form
\begin{equation} \label{eq:terms_R}
  c (1+b)^{1-l} \prod_{j=1}^{k'}
    \partial_x^{\alpha_j} \partial_\xi^{\beta_j} \Bsub
    \prod_{j=k'+1}^l \partial_x^{\alpha_j} \partial_\xi^{\beta_j} b
\end{equation}
modulo $S^{-\infty}$ and with $\sum_j |\alpha_j| + k'
= \sum_j |\beta_j| + k' \geq k$.
By (\ref{eq:assumption_b1}) and (\ref{eq:assumption_b2})
each of these terms is a bounded family of symbols in
$S^{\gamma - (\sum_j |\alpha_j| + k')(1-2\frac{\gamma}{L})}
(\R^n \times \R^n)$.
Set $Q^{(k)}$ equal to the $- \half (1+b)^{-1/2}$ times the sum of the
terms with $\sum_j |\alpha_j| + k' = k$, and make it
selfadjoint by averaging with the adjoint.
Then $Q^{(k)}$ is a sum of terms of the form
(cf.\ Lemma 18.1.7 on adjoints of pseudodifferential operators
in \cite{hormander85a})
\begin{equation}
  c (1+b)^{1/2-l} \prod_{j=1}^{k'}
    \partial_x^{\alpha_j} \partial_\xi^{\beta_j} \Bsub ,
    \prod_{j=k'+1}^l \partial_x^{\alpha_j} \partial_\xi^{\beta_j} b
\end{equation}
modulo $S^{-\infty}$.
It follows that $R^{(k+1)}$ is again a sum of terms
(\ref{eq:terms_R}), now with $\sum_j |\alpha_j| + k'
= \sum_j |\beta_j| + k' \geq k + 1$.
The $Q^{(k)}$ are bounded families of symbols in
$S^{\gamma/2 - k(1-2\frac{\gamma}{L})}(\R^n \times \R^n)_
{1-\frac{\gamma}{L},\frac{\gamma}{L}}$,
with $\partial_z^j Q^{(k)}$ in
$S^{\gamma/2 + \frac{j \gamma}{L} - k(1-2\frac{\gamma}{L})}_
{1-\frac{\gamma}{L},\frac{\gamma}{L}}(\R^n \times \R^n)$.
It follows that the asymptotic sum $Q$ exists and satisfies
$Q^2 = 1 + B + R$, with $R \in S^{-\infty}$.
\end{proof}

Next we consider the well-posedness of the Cauchy problem
\begin{align} \label{eq:Cauchy_problem_wellposedness}
  Pu = {}& f , \qquad 0 < z < Z ; &
   u = {}& u_0 , \qquad \text{ when } z = 0 .
\end{align}
For convenience we have set $z_0 = 0$ here.
To show the well-posedness of this Cauchy problem we use the previous lemma
and follow the argument in \cite{hormander85a}, Lemma 23.1.1
and Theorem 23.1.2, with minor modifications.
Let $\widetilde{\gamma} = \max(1,\gamma)$. We have the following
equivalents to Lemma 23.1.1 and Theorem 23.1.2.

\begin{lemma} \label{lem:energy_estimate_wellposedness}
Suppose $A$ and $B$ satisfy (i), (ii) and (iii) with $2\gamma < L$.
If $s \in \R$ and if $\lambda$ is larger than some number depending on
$s$, we have for every $u \in C^1([0,Z];H^s) \cap
C^0([0,Z];H^{s+\widetilde{\gamma}})$ and $p \in [1,\infty]$
\begin{equation} \label{eq:energy_estimate_wellposedness}
  \bigg( \half \int_0^Z \| \e^{-\lambda z} u(z,\cdot) \|_{H^s}^p \bigg)^{1/p}
  \leq \| u(0,\cdot) \|_{H^s}
    + 2 \int_0^Z \e^{-\lambda z} \| P u \|_{H^s} \, \dd z ,
\end{equation}
with the interpretation as a maximum when $p = \infty$.
\end{lemma}

\begin{theorem} \label{th:well-posedness}
Suppose $A$ and $B$ satisfy (i), (ii) and (iii) with $2\gamma < L$.
For every $f \in L^1(]0,Z[; H^s)$ and $u_0 \in H^s$, there is then a
unique solution $u \in C([0,T]; H^s)$ of the Cauchy problem
(\ref{eq:Cauchy_problem_wellposedness}), and
(\ref{eq:energy_estimate_wellposedness}) remains valid for this
solution.
\end{theorem}

\begin{proof}[Proof of Lemma~\ref{lem:energy_estimate_wellposedness} and
  Theorem~\ref{th:well-posedness}]
With the operators $Q$ and $R$ from Lemma~\ref{lem:root_B} we find
that for some constant $c$
\begin{align}
  \Re ((\ii A(z,x,D_x) + B(z,x,D_x)) v , v)
    = {}& \Re ((\ii A(z,x,D_x) + R - 1) v, v) + (Q v, Qv)
    \nonumber\\
  \geq {}& - c (v,v),
    \qquad v \in H^{\widetilde{\gamma}} ,
\end{align}
where the sharp G\aa{}rding inequality is used to estimate $\Re
(\ii A(z,x,D_x) v,v)$. The case $s = 0$ now follows from the arguments of
the proof of Lemma 23.1.1 in \cite{hormander85a}. For $s \neq 0$ we
set $E_s(D_x) = (1 + |D_x|^2)^{s/2}$.
It follows from the composition formula for symbols and the assumptions
(\ref{eq:assumption_b1}), (\ref{eq:assumption_b2}) that
$E_s(D_x) A(z,x,D_x) E_{-s}(D_x)$ and $E_s(D_x) B(z,x,D_x) E_{-s}(D_x)$
satisfy the same assumptions as $A$ and $B$. For $s \neq 0$ the estimate
(\ref{eq:energy_estimate_wellposedness}) now follows from the same
estimate with $s$ replaced by $0$, $P$ replaced by $E_s(D_x) P
E_{-s}(D_x)$, and $u$ replaced by $E_s(D_x) u$.

For the theorem we follow the proof of Theorem~23.1.2 in
\cite{hormander85a}. In the proof of the uniqueness in Theorem~23.1.2
we replace $s-1$ by $s - \widetilde{\gamma}$. For the existence we
find that the constructed solution
$u \in C^1([0,Z];H^{s-2\widetilde{\gamma}})$. Thus (on page 388) $s$
is replaced by $s+2\widetilde{\gamma}$ instead of $s+2$, and the
approximating solution $u_\nu$ should be in
$C^1([0,Z];H^{s+\widetilde{\gamma}})$. The theorem then follows.
\end{proof}

\section{Symbol estimates for $\exp(-I)$\label{sec:symbol_estimates}}

In this section we will establish symbol estimates for functions
of the form $K \exp(-I) = K(z,x,\xi) \exp(-I(z,x,\xi))$,
$z \in [z_0,Z]$.
Differentiation of $\exp(-I)$ brings out a derivative of $I$ that
is of order $(1+|\xi|)^\gamma$. We will improve on this by using
the property (iii). We will first estimate the derivatives
of $I$ by powers of $(1+|\xi|)$ and $I$. After that symbol estimates are
obtained using the fact that $(1+y)^\delta \exp(-y)$ is bounded for
$y \geq 0$. We study the case $A=0$, where the
bicharacteristics of $P_0$ are the straight lines
$(x,\xi) = \text{constant}$ and the integral $I$ reduces to
\begin{equation} \label{eq:define_I_integral_a=0}
  I(z,x,\xi)
  = \int_{z_0}^z b(z',x,\xi) \, \dd z' .
\end{equation}

\begin{lemma} \label{lem:estimate_b_I1}
There are constants $C_{\alpha,\beta}$ such that
\begin{equation} \label{eq:estimate_der_I_noz}
  \partial_x^\alpha \partial_\xi^\beta I
  \leq C_{\alpha,\beta} (1+|\xi|)^{-|\beta|+\frac{|\alpha|+|\beta|}{L}\gamma}
    (1+I)^{1-\frac{|\alpha|+|\beta|}{L}} , \quad
    z \in [z_0,Z] .
\end{equation}
\end{lemma}

\begin{proof}
This is true automatically for $|\alpha| + |\beta| \geq L$.
So suppose $|\alpha| + |\beta| < L$. By definition
\begin{equation}  \label{eq:derivative_I_noz}
  \partial_x^\alpha \partial_\xi^\beta I
  = \int_{z_0}^z \partial_x^\alpha \partial_\xi^\beta b(z',x,\xi)
    \, \dd z' , \quad |\xi| > 1 .
\end{equation}
The H\"older inequality implies that
\begin{align} \label{eq:inequality_integral_b}
  \int_{z_0}^z
      \big((1+|\xi|)^{-\gamma} & (1+b) \big)^{1-\frac{|\alpha|+|\beta|}{L}}
      \, \dd z' \nonumber\\
    \leq {}&
      \big\| \big( (1+|\xi|)^{-\gamma} (1+b) \big)^%
            {1-\frac{|\alpha|+|\beta|}{L}} \big\|_%
            {L^{1/(1-\frac{|\alpha|+|\beta|}{L})}}
      \big\| I_{[z_0,z]} \big\|_{L^{\frac{L}{|\alpha|+|\beta|}}}
      \nonumber \\
    = {}& (z-z_0)^{\frac{|\alpha|+|\beta|}{L}}
        \left( \frac{z-z_0 + I}{(1+ |\xi|)^\gamma}\right)^%
           {1-\frac{|\alpha|+|\beta|}{L}} ,
      \nonumber \\
    \leq {}& C (z-z_0)^{\frac{|\alpha|+|\beta|}{L}}
        \left( \frac{1 + I}{(1+ |\xi|)^\gamma}\right)^%
           {1-\frac{|\alpha|+|\beta|}{L}} ,
\end{align}
where $I_{[z_0,z]}$ is the indicator function and
the $L^{\frac{L}{|\alpha|+|\beta|}}$ norm is taken on a
$z$-interval in $\R$, for fixed $(x,\xi)$.
This inequality and (\ref{eq:assumption_b1}) imply the
estimate
\begin{equation} \label{eq:I_xxi_precise}
  \partial_x^\alpha \partial_\xi^\beta I
  \leq C_{\alpha,\beta}
    (1+|\xi|)^{-|\beta|+\frac{|\alpha|+|\beta|}{L}\gamma}
    (z-z_0)^{\frac{|\alpha|+|\beta|}{L}}
    (1+I)^{1-\frac{|\alpha|+|\beta|}{L}} , \quad
    |\alpha|+|\beta| < L ,
\end{equation}
from which (\ref{eq:estimate_der_I_noz}) follows.
\end{proof}

\begin{lemma} \label{lem:I_bound_symbols1}
Suppose there is a constant $\mu$ such that for each
$\alpha,\beta,j$, there are constants $C_{\alpha,\beta,j}$ and
$p_{\alpha,\beta,j}$ such that
\begin{multline} \label{eq:special_symbol_I_bound1}
  \partial_x^\alpha \partial_\xi^\beta \partial_z^j K(z,x,\xi)
  \leq C_{\alpha,\beta,j}
    (1+|\xi|)^{\mu-|\beta|+(\frac{|\alpha|+|\beta|}{L}+j)\gamma}
     (1+I(z,x,\xi))^{p_{\alpha,\beta,j}} , z \in [z_0,Z] .
\end{multline}
Then $K(z,x,\xi) \exp(-I(z,x,\xi))$ is a bounded family
of symbols in $S^{\mu}_{1-\frac{\gamma}{L},\frac{\gamma}{L}}
(\R^n \times \R^n)$, with $\partial_z^j (K \exp(-I))$
in $S^{\mu+j\gamma}_{1-\frac{\gamma}{L},\frac{\gamma}{L}}
(\R^n \times \R^n)$.
\end{lemma}

\begin{proof}
For the derivative of an exponential we have
\begin{multline} \label{eq:derivative_exponential}
  \partial_x^\alpha \partial_\xi^\beta \partial_z^j (K \exp(-I))
  =
  \sum_{m=1}^{|\alpha|+|\beta|+j}
    \sum_{\alpha^{(1)}+\ldots+\alpha^{(m)}=\alpha}
    \sum_{\beta^{(1)}+\ldots+\beta^{(m)}=\beta}
    \sum_{j^{(1)}+\ldots+j^{(m)}=j}\\
  c_{m,(\alpha^{(1)},\ldots,\alpha^{(m)}),(\beta^{(1)},\ldots,\beta^{(m)}),%
        (j^{(1)},\ldots,j^{(m)})} \exp(-I)\\
  \partial_x^{\alpha^{(1)}} \partial_\xi^{\beta^{(1)}}
      \partial_z^{j^{(1)}} K(z,x,\xi)
    \prod_{l=2}^m
      \partial_x^{\alpha^{(l)}} \partial_\xi^{\beta^{(l)}}
      \partial_z^{j^{(l)}} (-I) .
\end{multline}
Here $\sum_{j^{(1)}+\ldots+j^{(m)}=j}$ is the sum over all $m$-vectors with
non-negative integer components satisfying $j^{(1)}+\ldots+j^{(m)}=j$, and
for $\sum_{\alpha^{(1)}+\ldots+\alpha^{(m)}=\alpha}$ the $\alpha^{(l)}$ are
itself multi-indices, and the sum is such that
$(\alpha^{(l)},\beta^{(l)},j^{(l)}) \neq 0$ for $l\geq 2$.
From the inequality (\ref{eq:estimate_der_I_noz}) and the fact that
$\partial_x^\alpha \partial_\xi^\beta \partial_z^j I \leq
C (1+|\xi|)^{\gamma-|\beta|}$ it follows that
\begin{multline}
  \big| \partial_x^{\alpha^{(1)}} \partial_\xi^{\beta^{(1)}}
      \partial_z^{j^{(1)}} K(z,x,\xi)
    \prod_{l=2}^m
    \partial_x^{\alpha^{(l)}} \partial_\xi^{\beta^{(l)}}
    \partial_z^{j^{(l)}} (-I) \big| \\
  \leq
    C (1+|\xi|)^{\mu-|\beta|+(\frac{|\alpha|+|\beta|}{L}+j)\gamma} (1+I)^q .
\end{multline}
for some $q$. Since $(1+y)^\delta \exp(-y)$ is bounded for any $\delta$,
$y \geq 0$ it follows that
\begin{equation}
  \partial_x^\alpha \partial_\xi^\beta \partial_z^j (K \exp(-I))|
  \leq C (1+|\xi|)^{\mu-|\beta|+(\frac{|\alpha|+|\beta|}{L}+j)\gamma} .
\end{equation}
This completes the proof.
\end{proof}

The lemma shows in particular that the function $\exp(-I)$ is a
bounded family of
symbols in $S^{0}_{1-\frac{\gamma}{L},\frac{\gamma}{L}}
(\R^n \times \R^n)$, with $\partial_z^j \exp(-I)$ a bounded
family in $S^{j\gamma}_{1-\frac{\gamma}{L},\frac{\gamma}{L}}
(\R^n \times \R^n)$.

The class of $K$ satisfying the assumption of the lemma is closed
under multiplication and taking derivatives. That is, if $K'$ satisfies
the assumptions of Lemma~\ref{lem:I_bound_symbols1}, with constant $\mu'$,
then its derivative $\partial_x^\alpha \partial_\xi^\beta \partial_z^j
K'(z,x,\xi)$ satisfies the assumption with $\mu = \mu' - |\beta| +
(\frac{|\alpha|+|\beta|}{L} + j) \gamma$.
If $K''$ satisfies the assumption of the lemma with constant $\mu''$,
then the product $K' K''$ satisfies the assumptions with constant
$\mu' + \mu''$.

\section{Parametrix for the case $a=0$\label{sec:parametrix_A=0}}

In this section we prove Theorem~\ref{th:main_theorem} for the case
$A=0$, for which the operator $E_0$ is given by $E_0(z,z_0) =
\operatorname{Id}$.
An important part of the proof is an order by order construction.
We first prove a lemma that is important for the induction step.
As usual $\#$ denotes the composition of symbols.
From now on we will write $B^{(0)} = b, B^{(1)} = B - b$.

\begin{lemma} \label{lem:remainder_symbol_estimate1}
Assume $2\gamma < L$.
Suppose $K = K(z,x,\xi)$ is a smooth family of symbols in
$S^m(\R^n \times \R^n)$ that satisfies the
estimate (\ref{eq:special_symbol_I_bound1}) with constant $\mu$.
Then $B \# (K \exp(-I))$ can be written as an
asymptotic sum
\begin{equation} \label{eq:symbol_diff}
    B(z) \# \big( K(z) \exp(-I(z))
    = \sum_{j=0}^\infty M^{(j)} \exp(-I(z)) + R ,
\end{equation}
with $R$ a smooth family in $S^{-\infty}$, $M^{(0)} = b K$
and $M^{(j)}$ a smooth family of symbols in
$S^{m+\gamma+j(\gamma-1)}(\R^n \times \R^n)$ that satisfies
(\ref{eq:special_symbol_I_bound1}) with constant
$\mu' = \mu + \gamma - j + \frac{j}{L}\gamma$.
The integral $\int_{z_0}^z M^{(j)}(z') \, \dd z'$
satisfies (\ref{eq:special_symbol_I_bound1}) with
$\mu' = \mu - j(1 - \frac{2}{L}\gamma)$.
\end{lemma}

\begin{proof}
By the composition formula the symbol
$B(z) \# (K(z) \exp(-I(z)))$ is given by an asymptotic sum
(with $j=0$ or $1$)
\begin{align}
  \sum_{\alphatilde,j}
            \frac{(-\ii)^{|\alphatilde|}}{\alphatilde!}
    \partial_\xi^\alphatilde B^{(j)}(z)
      \partial_x^\alphatilde (K(z)\exp(-I(z)) + R,
\end{align}
with $R$ a smooth family in $S^{-\infty}(\R^n \times \R^n)$.
We let
\begin{equation} \label{eq:M_k_sum}
  M^{(k)}(z) = \sum_{\alphatilde,j,|\alphatilde| + j = k}
            \frac{(-\ii)^{|\alphatilde|}}{\alphatilde!}
        \exp(I(z))
    \partial_\xi^\alphatilde B^{(j)}(z)
        \partial_x^\alphatilde (K(z)\exp(-I(z)) ,
\end{equation}
so that (\ref{eq:symbol_diff}) is satisfied.
Each term in the sum (\ref{eq:M_k_sum}) is a product of a constant,
a factor $\partial_\xi^\alphatilde B^{(j)}(z)$, a factor
$\partial_x^{\alpha'} K(z)$ and a factor $\exp(I(z))
\partial_x^{\alpha''} \exp(-I(z))$, with
$\alpha'+\alpha''=\alphatilde$. These are smooth families of symbols
of order $\gamma-j-|\alphatilde|$, $\mu$ and
$|\alpha''| \gamma$, which shows the first statement about the
$M^{(j)}$. Also, they satisfy the
assumption of Lemma~\ref{lem:I_bound_symbols1} with constant
$\mu$ equal to $\gamma-|\alphatilde|-j$, $\mu +\frac{|\alpha'|}{L}\gamma$
and $\frac{|\alpha''|}{L}\gamma$, respectively.
The remarks following the proof of Lemma~\ref{lem:I_bound_symbols1}
about the multiplication of such functions show that each term
satisfies the assumptions of this lemma with constant
$\mu + \gamma - |\alphatilde| - j + \frac{|\alphatilde|}{L}\gamma$.
Since for $M^{(k)}$ we have $k = |\alphatilde|+j$ it follows that
the $M^{(k)}$ satisfies (\ref{eq:special_symbol_I_bound1}) with constant
$\mu' = \mu + \gamma - j + \frac{j}{L}\gamma$.

Let $K'$ be given by $K' = \int_{z_0}^z M^{(l)}(z',x,\xi) \, \dd
z'$, for some $l$. For $K'$ we must estimate the multiple derivative
$\partial_x^\alpha \partial_\xi^\beta \partial_z^j K'(z,x,\xi)$
If $j \neq 0$, then this is equal to a multiple
derivative of $M^{(l)}$, given by $\partial_x^\alpha \partial_\xi^\beta
  \partial_z^{j-1} M^{(l)}$, and the result follows from
the result already proven for $M^{(l)}$.
Next suppose that $j = 0$. In this case
$\partial_x^\alpha \partial_\xi^\beta \partial_z^j K'(z,x,\xi)$
is a sum of terms
\begin{equation}
  c \int_{z_0}^z \partial_x^{\alpha'+\alphatilde} \partial_\xi^{\beta'} B^{(j)}
    \exp(I) \partial_x^{\alpha''} \partial_\xi^{\beta''+\alphatilde}
    (K\exp(-I)) \, \dd z' ,
\end{equation}
where $c$ is a constant and $\alpha'+\alpha''=\alpha, \beta'+\beta'' =
\beta$. For such term we can put outside the integral a factor
$C (1+|\xi|)^{\mu - |\beta''|-|\alphatilde| +
\frac{|\alpha''|+|\beta''|+|\alphatilde|}{L}\gamma} (1+I)^p$
that is an upperbound for $\exp(I) \partial_x^{\alpha''}
\partial_\xi^{\beta''+\alphatilde} (K \exp(-I))$.
Thus we obtain
\begin{multline} \label{eq:integral_parametrix_lemma1}
  \partial_x^\alpha \partial_\xi^\beta K' \leq
    C (1+|\xi|)^{\mu - |\beta''|-|\alphatilde|
    + \frac{|\alpha''|+|\beta''|+|\alphatilde|}{L}\gamma}
    (1+I)^p \\
  \times (1+|\xi|)^{\gamma-|\beta'|-j}
    \int_{z_0}^ z \big((1+|\xi|)^{-\gamma} (1+b)\big)^%
            {1 - \frac{|\alpha'+\alphatilde|+|\beta'| +2j}{L}} \, \dd z',
\end{multline}
if $|\alpha'|+|\alphatilde|+|\beta'| +2j < L$.
The integral on the right hand side can be estimated as in
(\ref{eq:inequality_integral_b}).  Since $l = j + |\alphatilde|$
this yields that $K'$ satisfies (\ref{eq:special_symbol_I_bound1})
with $\mu$ given by $\mu' = \mu - l(1 - \frac{2}{L}\gamma)$.
If $|\alpha'+\alphatilde|+|\beta'| +2j \geq L$, then
\begin{equation} \label{eq:integral_parametrix_lemma2}
  \partial_x^\alpha \partial_\xi^\beta K' \leq
    C (1+|\xi|)^{\mu - |\beta''|-|\alphatilde|
    + \frac{|\alpha''|+|\beta''|+|\alphatilde|}{L}\gamma + \gamma-|\beta'|-j}
    (1+I)^p .
\end{equation}
Then $(1+|\xi|)^\gamma \leq
(1+|\xi|)^{\frac{|\alpha'|+|\alphatilde|+|\beta'|+2j}{L}\gamma}$, so
that again $K'$ satisfies (\ref{eq:special_symbol_I_bound1})
with $\mu$ given by $\mu' = \mu - l(1 - \frac{2}{L}\gamma)$.
\end{proof}

\begin{proof}[Proof of Theorem~\ref{th:main_theorem} for the case $A=0$]
% cf.\ \cite{hormander85a}, th 18.1.35

In this case $E_0(z,z_0) = \operatorname{Id}$, and we can write $E =
W$, where we view the family of pseudodifferential operators $W$ as
mapping functions of $x \in \R^n$ to functions of
$(z,x) \in [z_0,Z] \times \R^n$. We construct a bounded family
of pseudodifferential operators $\What(z)$, $z \in [z_0,Z]$ in
$S^0_{1-\frac{\gamma}{L},\frac{\gamma}{L}}$, with $\partial_z^j
\What(z)$ a bounded
family in $S^{j\gamma}_{1-\frac{\gamma}{L},\frac{\gamma}{L}}$ satisfying
\begin{equation} \label{eq:Winit}
  \What(z_0,x,D_x) = \operatorname{Id} ,
\end{equation}
and
\begin{equation} \label{eq:PW}
  P \What \text{ is a bounded map }
    H^s(\R^n) \rightarrow C^k([z_0,z_0+Z],H^{s+l}(\R^n)) ,
\end{equation}
for any $k,s,l$.
After that we show that $W - \What$ is a smoothing operator.
From the construction of $\What$ it will also follow that $W$ has
the properties described in the theorem.

The operator $\What$ will be constructed as an asymptotic sum of
$W^{(j)}$ with decreasing order $\What = \sum_{j=0}^\infty W^{(j)}$.
The symbols $W^{(j)}(z,x,\xi)$ will be of the form
\begin{equation} \label{eq:form_W_j}
  W^{(j)}(z,x,\xi) = K^{(j)}(z,x,\xi) \exp(-I(z,x,\xi)) ,
\end{equation}
with $K^{(0)} = 1$ and the $K^{(j)}$, $j \geq 1$ to be determined. We
will assume that the $K^{(j)}$ will satisfy the assumptions of
Lemma~\ref{lem:I_bound_symbols1} for
$\mu = -j(1-\frac{2}{L}\gamma)$, and are symbols of order $j (2\gamma - 1)$.

The compositions $\partial_z \What$ and $B \What$ are again a family of
pseudodifferential operators. We have
\begin{equation}
  \partial_z W^{(j)} = {} \Op \big(
    \pdpd{K^{(j)}}{z} \exp(-I) - b K^{(j)} \exp(-I) \big) .
\end{equation}
Hence
\begin{align} \label{eq:P_Wj}
  P W^{(j)} = {}&\Op (\pdpd{K^{(j)}}{z} \exp(-I)) \nonumber\\
  & + \Op(B(z,x,\xi)) \Op (K^{(j)} \exp(-I))
    - \Op(b(z,x,\xi) K^{(j)} \exp(-I)) .
\end{align}
We define operators $R = P\What$ and $R^{(k)}$, given by
\begin{equation}
  R^{(k)} = P \sum_{j=0}^{k-1} W^{(j)} .
\end{equation}
We denote by $M^{(j,k)}$ the operators $M^{(k)}$ from the previous
lemma, applied to $K = K^{(j-1)}$.
For $k = 1$ we can write, using (\ref{eq:P_Wj})
\begin{equation} \label{eq:R_k_series}
  R^{(k)} = \sum_{l=0}^\infty r^{(k,l)} \exp(-I) ,
\end{equation}
with
\begin{align}
  r^{(1,l)} = {}& M^{(1,l)} , \quad l \geq 1 , \\
  r^{(1,0)} = {}& M^{(1,0)} + \pdpd{K^{(0)}}{z} - b = 0 .
\end{align}
For $k > 1$ we have (\ref{eq:R_k_series}) if we set
\begin{align}
  r^{(k,k-1+l)} = {}& r^{(k-1,k-1+l)} + M^{(k,l)} , \\
  r^{(k,k-1)} = {}& r^{(k-1,k-1)} + M^{(k,0)}
      + \pdpd{K^{(k-1)}}{z} - b K^{(k-1)} = r^{(k-1,k-1)} +
      \pdpd{K^{(k-1)}}{z}.
      \label{eq:rkk-1}
\end{align}

Assume that $r^{(k,j)} = 0$, $j=0,\ldots, k-1$. Then let
\begin{equation}
  K^{(k)}(z) = - \int_{z_0}^z r^{(k,k)}(z') \, \dd z' .
\end{equation}
By the previous lemma this is a bounded family of symbols in
$S^{k(2\gamma-1)}$. It also satisfies the assumptions of
Lemma~\ref{lem:I_bound_symbols1} with
$\mu = -k(1-\frac{2}{L}\gamma)$. By (\ref{eq:rkk-1}) it follows
that then $r^{(k+1,k)} = 0$. Thus by induction we find a series of
$K^{(k)}$ such that $R^{(k)}$ is a bounded family of symbols in
$S^{\gamma-\frac{\gamma}{L}-k(1-\frac{2}{L}\gamma)}_%
{1-\frac{\gamma}{L},\frac{\gamma}{L}}$,
for $z \in [z_0,Z]$. It follows that, for fixed $z$, $R$ is continuous
$H^s(\R^n)$ to $H^{s+l}(\R^{n+1})$, uniformly in $z$.
The operators $\partial_z^j R$ are bounded families of symbols in
$S^{(j+1)\gamma-\frac{\gamma}{L}}(\R \times \R)$.
The terms in their asymptotic expansions also vanish.
Hence is $\partial_z^j R$ is also continuous $H^s(\R^n)$ to
$H^{s+l}(\R^n)$, for any $s,l$.
This shows (\ref{eq:PW}).

By definition $W^{(0)}(z_0) = 1$ and $W^{(j)} = 0$, $j \geq 1$. We can assume
the asymptotic sum of symbols is such that (\ref{eq:Winit}) is satisfied.
(Because $\What$ is a bounded family, with continuous symbol, it follows that
$W(z)u_0$ is a continuous function of $z$ with values in $H^s$,
for any $u_0 \in H^s$.)
%Let $\psi : \R^n \rightarrow [0,1]$ have bounded support
%and be equal to $1$ on a neighborhood of zero. We define
%$\psi_\eta(\xi) = \psi(\eta \xi)$.
%Let $\epsilon > 0$. For $\eta>0$ sufficiently small we have that
%\begin{equation}
%  \| ( 1 - \psi_\eta) u_0 \|_{H^s} < \frac{\epsilon}{4 \min(C_{s,0},1)} ,
%\end{equation}
%with $C_{s,0}$ as in (\ref{eq:Wj_uniformbound}).
%The symbol $W^{(0)}(z,z_0,x,\xi) \psi_\eta(\xi)$ is uniformly
%continuous as a function of $z,z_0$ with values in
%$S^0_{1-\frac{\gamma}{L},\frac{\gamma}{L}}$. It follows that
%for $W^{(0)}(z,z_0) \psi_\eta u_0 \rightarrow u_0$ in $H^s(\R^n)$
%for $u_0 \in H^s(\R^n)$. Let $\delta$ be such that
%$\|W^{(0)}(z,z_0) \psi_\eta u_0 - \psi_\eta u_0\|_{H^s(\R^n)} \leq
%\frac{\epsilon}{4}$ for $z < z_0 + \delta$. Then for $z < z_0 +
%\delta$ we have
%\begin{align}
%  \| EW{(0)}(z,z_0) u_0 - u_0 \|_{H^s(\R^n)}
%  \leq {}& \| W^{(0)}(z,z_0) (1-\psi_\eta) u_0 - (1-\psi_\eta) u_0 \|_{H^s(\R)}
%      \nonumber\\
%      + \| W^{(0)}(z,z_0) \psi_\eta u_0 - \psi_\eta u_0 \|_{H^s(\R^n)}
%      \nonumber\\
%  \leq {}& \frac{3\epsilon}{4} .
%\end{align}

Equation (\ref{eq:PW}) and the energy estimates for this problem
(see section~\ref{sec:wellposedness}) imply that
\begin{equation} \label{eq:estimate_difference_E_W}
  W - \What \in S^{-\infty}(([z_0,Z] \times \R^n ) \times \R^n)
\end{equation}
It follows that the equations (\ref{eq:parametrix_left_and_right})
are satisfied.

From the construction of $\What$ it follows that (\ref{eq:Wprop1}) and
(\ref{eq:Wprop2}) are satisfied. It is clear that when $I(z,x,\xi) =
0$, then $b(z',x,\xi) = 0$, $z' \in [z_0,z]$, and then all the
$K^{(j)}(z,x,\xi)$ are zero.
\end{proof}

\section{The case $A \neq 0$\label{sec:complete_result}}

In this section we complete the proof Theorem~\ref{th:main_theorem},
by proving the case $A \neq 0$. The general case will be derived from
the case $A = 0$, by using Egorov's theorem (see e.g.\
\cite{taylor81}, p.\ 147). Consider the transformed function
$\widetilde{u}$ defined from $u$ by
\begin{equation}
  \widetilde{u}(z,\cdot) = E_0(z,z_0)^{-1} u(z) .
\end{equation}
Of course we have
\begin{equation}
  \pdpd{E_0}{z}(z,z_0) = \ii A(z,x,D_x) E_0(z,z_0) .
\end{equation}
We will use the notation
\begin{equation} \label{eq:conjugate_B}
  \widetilde{B}(z,x,D_x)
  = E_0(z,z_0)^{-1} B(z,x,D_x) E_0(z,z_0) .
\end{equation}
The pseudodifferential equation now becomes
\begin{equation} \label{eq:tilde_PsDE}
  \big( \partial_z + \widetilde{B}(z,x,D_x) \big)
  \widetilde{u} = 0 .
\end{equation}
By Egorov's theorem $\widetilde{B}$ is a pseudodifferential operator of
order $\gamma$ with homogeneous real non-negative principal symbol
\begin{equation} \label{eq:pullback_b}
  \widetilde{b} = b \circ \Phi_{z,z_0} .
\end{equation}
We have
\begin{equation}
  \widetilde{I}(z,x,\xi)
  = \int_{z_0}^z \widetilde{b}(z',x,\xi) \, \dd z' .
\end{equation}

The following lemma states that the properties (\ref{eq:assumption_b1}),
(\ref{eq:assumption_b2}) and (\ref{eq:form_of_W_left}),
(\ref{eq:Wprop1}), (\ref{eq:Wprop2}) are conserved under the mapping
$B \mapsto \widetilde{B}$.

\begin{lemma} \label{lem:FIO_conjugate}
The symbol $\widetilde{B}$, defined by
(\ref{eq:conjugate_B}) is a
smooth family of symbols in $S^\gamma(\R^n \times \R^n)$ satisfying
(\ref{eq:assumption_b1}) and (\ref{eq:assumption_b2}) with
$\widetilde{b}$ instead of $b$, if and only if $B$ has these properties.
Let $W(z,x,D)$ be a bounded family of pseudodifferential
operators. Then $\widetilde{W}(z) = \exp(z,z_0)^{-1} W(z)
\linebreak[2]\exp(z,z_0)$
satisfies (\ref{eq:form_of_W_left}), (\ref{eq:Wprop1}),
(\ref{eq:Wprop2}) if and only if $W$ satisfies
(\ref{eq:form_of_W_left}), (\ref{eq:Wprop1}), (\ref{eq:Wprop2}).
\end{lemma}

\begin{proof}
We use Egorov's theorem in the form given by Taylor \cite{taylor81}, p.\
147. Taylor assumes that the subprincipal symbol $A - a$ is
polyhomogeneous, but it can be checked from the proof that this
assumption may be omitted, and that the result applies to our case as
well. The map $Q \mapsto \widetilde{Q} = E(z,z_0)^{-1} Q E(z,z_0)$ maps a
bounded set of symbols $Q(x,\xi)$ to a bounded set of symbols, where
the principal symbol is given by $\widetilde{q} = q \circ \Phi_{z,z_0} $.

To apply this to a family of symbols $B = B(z,x,\xi)$, observe that
the derivatives $\partial_z^j \widetilde{B}$ are given by
\begin{equation} \label{eq:z_der_Btilde}
  \partial_z^j \widetilde{B}
  = E(z,z_0)^{-1} \big(\partial_z - \ii [A, \cdot ]\big)^j B E(z,z_0) .
\end{equation}
It follows that if $B$ is a smooth family of symbols in
$S^\gamma(\R^n \times \R^n)$ then $\widetilde{B}$ is a
smooth family of symbols in $S^\gamma(\R^n \times \R^n)$.

It follows from (\ref{eq:pullback_b}) that $\widetilde{B}$ has
homogeneous, real, non-negative principal symbol $\widetilde{b}$.
To establish that $\widetilde{B}, \widetilde{b}$ satisfy the properties
(\ref{eq:assumption_b1}), (\ref{eq:assumption_b2}) we recall the
construction of the asymptotic series for $\widetilde{B}$ in the proof of
Egorov's theorem in \cite{taylor81}, p.\ 147. Consider the transformation
of a symbol $Q$ independent of $z$. From the transformation
$\widetilde{Q} = E(z,z_0)^{-1} Q E(z,z_0)$ a $z_0$-family of
operators $\widetilde{Q}$ is obtained, that satisfies the differential
equation
\begin{equation}
  \partial_{z_0} \widetilde{Q}(z_0)
  = \ii [ A(z_0,x,D), \widetilde{Q}(z_0)] ,\quad \widetilde{Q}(z) = Q .
\end{equation}
\newcommand{\Qtilde}{\widetilde{Q}}
The asymptotic series is obtained by solving a series of differential
equations for the $\Qtilde^{(j)}$, $j=0,1,\ldots$
\begin{equation} \label{eq:diff_eq_q0}
  (\partial_{z_0} - H_a) \Qtilde^{(j)}(z_0,x,\xi) = a_{j-1} ,
\end{equation}
with initial condition
\begin{equation}
  \widetilde{Q}^{(j)}(z_0,x,\xi) =
  \left\{ \begin{array}{ll} Q(x,\xi) , \text{ at } z_0 = z
                                       & j=0\\
                            0 ,        \text{ at } z_0 = z
                                       & j\geq 1 ,\end{array} \right.
\end{equation}
where $a_{-1} = 0$ and
\begin{equation}
  a_j = \{ A - a, \Qtilde^{(j)} \}
      + \sum_{|\alpha|\geq2} \frac{\ii^{-|\alpha|+1}}{\alpha!}
        (\partial_\xi^\alpha A \partial_x^\alpha \Qtilde^{(j)} -
        \partial_\xi^\alpha \Qtilde^{(j)} \partial_x^\alpha A) .
\end{equation}
Thus $\Qtilde^{(0)}$ is constant along the integral curves of
$\partial_z - H_a$, hence we have
\begin{equation}
  \Qtilde^{(0)}(z_0,x,\xi) = Q(\Phi_{z_0,z}^{-1}(x,\xi)) .
\end{equation}
For $j \geq 1$ we find
\begin{equation}
  \Qtilde^{(j)}(z_0,x,\xi)
  = - \int_{z_0}^z a_{j-1}(z',\Phi(z_0,z')^{-1}(x,\xi)) .
\end{equation}
It follows that we can write the asymptotic series for $Q$ as
\begin{equation} \label{eq:transformed_series_Q}
  \bigg( \sum_{\alpha,\beta}
    S_{\alpha,\beta} \partial_x^\alpha \partial_\xi^\beta Q \bigg)
    \circ \Phi_{z_0,z}^{-1} ,
\end{equation}
where $S_{\alpha,\beta}$ are smooth families of symbols of order
$-|\alpha| + \lfloor \frac{|\alpha|+|\beta|}{2} \rfloor$ and
$S_{0,0} = 1$.

This shows that (\ref{eq:assumption_b1}), (\ref{eq:assumption_b2})
are satisfied at least when the $z$-component of $\alpha$ is $0$.
For terms with a non-zero number of $z$-derivatives it follows
using (\ref{eq:z_der_Btilde}).

This shows the if part. The map $B \mapsto \widetilde{B}$ has an
inverse given by $B = E(z,z_0) \widetilde{B}
\linebreak[2] E(z,z_0)^{-1}$, which
satisfies a similar differential equation in $z$. The only
if part is obtained by applying similar arguments.

The statement about $W$ follows from (\ref{eq:form_of_W_left})
and the expression for the asymptotic series
(\ref{eq:transformed_series_Q}).
% define the $\widetilde{K}^{(j)}$ ?????
\end{proof}

\begin{proof}[Proof of Theorem~\ref{th:main_theorem} for the general case]
By the previous lemma, and the proof of the theorem for $A = 0$ it
follows that (\ref{eq:tilde_PsDE}) has a solution operator
$\widetilde{E}$ that is a family of pseudodifferential operators
$\widetilde{W}(z,x,D)$ with the properties of the theorem.
By the definition of $\widetilde{u}$ it follows that
\begin{equation}
  E(z,z_0) = E_0(z,z_0) \widetilde{E}(z,z_0) .
\end{equation}
This yields the second equality in (\ref{eq:parametrix_left_and_right}).
Obviously we have
\begin{equation}
  E = E_0 \widetilde{E} E_0^{-1} E_0 .
\end{equation}
By Lemma~\ref{lem:FIO_conjugate} is a family of pseudodifferential
operators $W(z,x,D)$ with the properties of given in the theorem.
This yields the first equality in (\ref{eq:parametrix_left_and_right}).
\end{proof}

\section{The Fourier integral operator property\label{sec:FIO_property}}

In this section wel establish the Fourier integral operator property.
So far we had that for each $z$ derivative the bounds for symbols of
the form $K \exp(-I)$ increased with a
factor $(1+|\xi|)^\gamma$, uniformly in $z$. To obtain symbols
estimates also w.r.t.\ the $z$-derivatives we establish improved bounds in
the following lemmas.

\begin{lemma} \label{lem:estimate_b_I2}
Let $b$ be as in (ii), (iii). Suppose $I$ is given by
(\ref{eq:define_I_integral_a=0}). Then there are constants $C,C'$ such that
\begin{equation} \label{eq:estimate_b_I}
  b(z',x,\xi)
  \leq C (z-z_0)^{-\frac{L}{L+1}} (1+|\xi|)^{\frac{1}{L+1}\gamma}
    I(z,x,\xi)^{\frac{L}{L+1}}
\end{equation}
and
\begin{equation} \label{eq:b_joint_estimate}
  b(z',x,\xi) \leq C' (1+|\xi|)^\gamma
    (1+(z-z_0)(1+\xi)^\gamma)^{-1+\frac{1}{L+1}}
    (1+I(z,x,\xi))^{1-\frac{1}{L+1}} .
\end{equation}
when $z\in ]z_0,Z]$ and $z' \in [z_0,z]$.
There are constants $C_{\alpha,\beta,j}$ such that
\begin{equation} \label{eq:estimate_der_I_derz}
  \partial_x^\alpha \partial_\xi^\beta \partial_z^j I
  \leq C_{\alpha,\beta,j}
    (1+|\xi|)^{j \gamma-|\beta|}
    (1+(z-z_0)(1+|\xi|)^\gamma)^{-j+\frac{j+|\alpha|+|\beta|}{L}}
    (1+I)^{1-\frac{j+|\alpha|+|\beta|}{L}} ,
\end{equation}
when $z\in ]z_0,Z]$.
\end{lemma}

\begin{proof}
An assumption on $b$ is that
\begin{equation}
  (1+|\xi|)^{-\gamma} | \pdpd{b}{z} |
  \leq C \big( (1+|\xi|)^{-\gamma}  b\big)^{1-\frac{1}{L}} .
\end{equation}
It follows from this that for some constant $C_1$
\begin{equation}
  ((1+|\xi|)^{-\gamma} b(z'))^{\frac{1}{L}} -
  ((1+|\xi|)^{-\gamma} b(z''))^{\frac{1}{L}} \leq C_1 | z' - z''| .
\end{equation}
Assuming for the moment $z' > z''$, it follows that
\begin{equation}
  ((1+|\xi|)^{-\gamma} b(z'')) \geq
    \left( ((1+|\xi|)^{-\gamma} b(z'))^{\frac{1}{L}} - C_1 (z'-z'')\right)^L .
\end{equation}
Let $\widetilde{z}_\pm = z' \pm \frac{((1+|\xi|)^{-\gamma}
b(z'))^{\frac{1}{L}}}{C_1}$. Denote $I(z,z_0,x,\xi) = \int_{z_0}^z
b(z',x,\xi) \, \dd z'$. By integrating the previous inequality over
$z''$ it follows that
\begin{equation} \label{eq:integrated_inequality}
  (1+|\xi|)^{-\gamma} I(z',\widetilde{z}_-,x,\xi)
  \geq \frac{1}{L+1}((1+|\xi|)^{-\gamma} b(z',x,\xi))^{\frac{L+1}{L}} .
\end{equation}
With $z'' > z'$ we find the same inequality for
$(1+|\xi|)^{-\gamma} I(\widetilde{z}_+,z',x,\xi)$ by a similar
argument.
Hence if $z-z_0 \leq C ((1+|\xi|)^{-\gamma} b(z',x,\xi))^{\frac{1}{L}}$, then
\begin{equation}
  \frac{I}{z-z_0} \geq \frac{I(\widetilde{z}_-,z')}{z'-\widetilde{z}_-}
  \geq C b ,
\end{equation}
while if $Z-z_0 \geq z-z_0
\geq C ((1+|\xi|)^{-\gamma} b(z',x,\xi))^{\frac{1}{L}}$, then by
(\ref{eq:integrated_inequality})
\begin{equation}
  (1+|\xi|)^{-\gamma} I \geq C ((1+|\xi|)^{-\gamma} b)^{\frac{L+1}{L}} .
\end{equation}
It follows that
\begin{equation}
  (1+|\xi|)^{-\gamma} b
  \leq C \left( \frac{I}{(z-z_0) (1+|\xi|)^\gamma} \right)^{\frac{L}{L+1}} .
\end{equation}
The inequality (\ref{eq:estimate_b_I}) follows from this.

Of course we also have that $b \leq C(1+|\xi|)^\gamma$, hence
\begin{equation}
  b \leq C (1+|\xi|)^{\frac{\gamma}{L+1}} (1+I)^{\frac{L}{L+1}}
    \min\big( (1+|\xi|)^\gamma, (z-z_0)^{-1} \big)^{\frac{L}{L+1}} .
\end{equation}
The minimum on the right hand side can be estimated by
\begin{equation}
  \min\big( (1+|\xi|)^\gamma, (z-z_0)^{-1} \big)
  \leq 2 \big( \frac{1}{(1+|\xi|)^\gamma} + z-z_0 \big)^{-1}
    = \frac{2 (1+|\xi|)^\gamma}{1+(1+|\xi|)^\gamma(z-z_0)} .
\end{equation}
Thus we have (\ref{eq:b_joint_estimate}).

Suppose first that $j+|\alpha|+|\beta| < L$.
If $j = 0$ the inequality (\ref{eq:estimate_der_I_derz}) follows from
(\ref{eq:I_xxi_precise}). In case $j\geq 1$ we have for the
derivatives w.r.t.\ $(z,x,\xi)$ of $I$
\begin{equation}     \label{eq:derivative_I_derz}
  \partial_x^\alpha \partial_\xi^\beta \partial_z^j I
  = \partial_x^\alpha \partial_\xi^\beta \partial_z^{j-1} b(z,x,\xi) .
\end{equation}
From (\ref{eq:assumption_b1}), (\ref{eq:b_joint_estimate}) and
(\ref{eq:derivative_I_derz}) it follows that
\begin{equation} \label{eq:estimate_der_I_derz_jprecise}
  \partial_x^\alpha \partial_\xi^\beta \partial_z^j I
  \leq C
    (1+|\xi|)^{\gamma-|\beta|}
    (1+(z-z_0)(1+|\xi|)^\gamma)^{\frac{-L-1+j+|\alpha|+|\beta|}{L+1}}
    (1+I)^{\frac{L+1-j+|\alpha|+|\beta|}{L+1}} .
\end{equation}
The inequality (\ref{eq:estimate_der_I_derz}) follows, using
that $\frac{1+I}{1+(z-z_0)(1+|\xi|)^{\gamma} }$ is bounded.
For any $(\alpha,\beta,j)$ there is a constant $C$ such that
\begin{equation}
  \partial_x^\alpha \partial_\xi^\beta \partial_z^j I
  \leq C (z-z_0)(1+|\xi|)^{\gamma-|\beta|} .
\end{equation}
This implies that (\ref{eq:estimate_der_I_derz}) is also true
when $j+|\alpha|+|\beta| \geq L$.
\end{proof}

Analogously as in Lemma~\ref{lem:I_bound_symbols1}, symbols estimates
for functions of the form $K \exp(-I)$ in
$S^\mu_{1-\frac{\gamma}{L},\frac{\gamma}{L}} ((]z_0,Z[ \times \R^n)
\times \R^n)$ follow if $K$ satisfies certain estimates in terms of
powers of $(1+I)$:

\begin{lemma} \label{lem:I_bound_symbols2}
Suppose $I$ is given by (\ref{eq:define_I_integral_a=0}).
Suppose that there are $\mu$ and $\kappa$ such that for each
$\alpha,\beta,j$ there are constants $C_{\alpha,\beta,j}$ and
$p_{\alpha,\beta,j}$ such that
\begin{multline} \label{eq:special_symbol_I_bound2}
  \partial_x^\alpha \partial_\xi^\beta \partial_z^j K(z,x,\xi)
  \leq C_{\alpha,\beta,j}
    (1+|\xi|)^{\mu+j\gamma-|\beta|}\\
    \times \big(1+(z-z_0)(1+|\xi|)^\gamma\big)^
        {\kappa-j+\frac{j+|\alpha|+|\beta|}{L}}
    (1+I(z,x,\xi))^{p_{\alpha,\beta,j}} , \quad
    z \in ]z_0,Z] .
\end{multline}
Then there are constants $C'_{\alpha,\beta,j}$ such that
\begin{multline}
  \partial_x^\alpha \partial_\xi^\beta \partial_z^j
    ( K(z,x,\xi) \exp(-I(z,x,\xi)))\\
  \leq C'_{\alpha,\beta,j}
    (1+|\xi|)^{\mu+j\gamma-|\beta|}
    \big(1+(z-z_0)(1+|\xi|)^\gamma\big)^
        {\kappa-j+\frac{j+|\alpha|+|\beta|}{L}}, \quad
    z \in ]z_0,Z] ,
\end{multline}
and hence $K\exp(-I)) \in S^{\mu+\gamma\kappa}
_{1-\frac{\gamma}{L},\frac{\gamma}{L}}((]z_0,Z[\times\R^n)\times\R^n)$.
\end{lemma}

The class of $K$ satisfying the assumption of the lemma for
some $\mu$ and $\kappa$ is closed under multiplication and taking
derivatives. That is, if $K'$ satisfies the assumptions of
Lemma~\ref{lem:I_bound_symbols2}, with constant $\mu'$ and $\kappa'$,
then its derivative $\partial_x^\alpha \partial_\xi^\beta \partial_z^j
  K'(z,x,\xi)$ satisfies the assumption with
$\mu = \mu' - |\beta|$ and
$\kappa = \kappa' - j + \frac{j+|\alpha|+|\beta}{L}$.
If $K''$ satisfies the assumption of the lemma with constant $\mu''$
and $\kappa''$, then the product $K' K''$ satisfies the assumptions
with constants $\mu' + \mu''$ and $\kappa' + \kappa''$.
From Lemma~\ref{lem:estimate_b_I2} it follows that $I$ satisfies the
assumption of Lemma~\ref{lem:I_bound_symbols2} with $\mu = \kappa =
0$, and $b$ with $\mu=\gamma$ and $\kappa=-1 + \frac{1}{L}$.

\begin{proof}[Proof of Theorem~\ref{th:FIO_property}]
To show that the $W^{(j)}$ are in $S^{-j(1-\frac{2\gamma}{L})}_%
{1-\frac{\gamma}{L},\frac{\gamma}{L}}((]z_0,Z[\times\R^n)\times\R^n)$,
we first give an improvement of Lemma~\ref{lem:remainder_symbol_estimate1}.
We will show that in fact the $M^{(j)}$ satisfy
(\ref{eq:special_symbol_I_bound2}) with constants
$(\mu+\gamma-j, -1 + \frac{1+j}{L})$,
and that the integral $\int_{z_0}^z M^{(j)}(z',z_0,x,\xi) \, \dd z'$
satisfies (\ref{eq:special_symbol_I_bound2}) with constants
$(\mu, \kappa + \frac{2j}{L})$.

Recall that the $M^{(j)}$ were defined in (\ref{eq:M_k_sum}).
Each term in the sum (\ref{eq:M_k_sum}) is a product of a constant,
a factor $\partial_\xi^\alphatilde B^{(j)}(z)$, a factor
$\partial_x^{\alpha'} K(z)$ and a factor $\exp(I(z,z_0))
\partial_x^{\alpha''} \exp(-I(z,z_0))$. These satisfy the
assumptions of Lemma~\ref{lem:I_bound_symbols2} with constant
$(\mu,\kappa)$ equal to
$(\gamma-|\alphatilde|-j, -1 + \frac{1+2j+|\alphatilde|}{L})$,
$(\mu, \kappa + \frac{|\alpha'|}{L})$ and
$(0,\frac{|\alpha''|}{L})$, respectively.
The remarks following the proof of Lemma~\ref{lem:I_bound_symbols2}
about the multiplication of such functions show that each term
satisfies the assumptions of this lemma with constants
$(\mu+\gamma-|\alphatilde|-j, -1 + \frac{1+2j+2|\alphatilde|}{L})$.
Since for $M^{(k)}$ we have $k = |\alphatilde|+j$ it follows that
the $M^{(k)}$ satisfies the
assumptions of Lemma~\ref{lem:I_bound_symbols2} with constants
$(\mu+\gamma-k, -1 + \frac{1+k}{L})$.

Let $K'$ be given by $K' = \int_{z_0}^z M^{(l)}(z',x,\xi) \, \dd
z'$, for some $l$. For $K'$ we must estimate the multiple derivative
$\partial_x^\alpha \partial_\xi^\beta \partial_z^j K'(z,x,\xi)$.
If $j \neq 0$, then this is equal to a multiple
derivative of $M^{(l)}$, given by $\partial_x^\alpha \partial_\xi^\beta
  \partial_z^{j-1} M^{(l)}$, and the result follows from
the result already proven for $M^{(l)}$. Next suppose that $j = 0$.
In this case
$\partial_x^\alpha \partial_\xi^\beta \partial_z^j K'(z,x,\xi)$
is a sum of terms
\begin{equation}
  c \int_{z_0}^z \partial_x^{\alpha'+\alphatilde} \partial_\xi^{\beta'} B^{(j)}
    \exp(I) \partial_x^{\alpha''} \partial_\xi^{\beta''+\alphatilde}
    (K\exp(-I)) \, \dd z',
\end{equation}
where $c$ is a constant and $\alpha'+\alpha''=\alpha, \beta'+\beta'' =
\beta$. For such term we can put outside the integral a factor
$C (1+|\xi|)^\mu (1+(z-z_0)(1+|\xi|)^\gamma)^{\kappa+
\frac{|\alpha''|+|\beta''|+|\alphatilde|}{L}} (1+I)^p$
that is an upperbound for $\exp(I) \partial_x^{\alpha''}
\partial_\xi^{\beta''+\alphatilde} (K \exp(-I))$.
Thus we obtain
\begin{multline} \label{eq:integral_parametrix_lemma}
  \partial_x^\alpha \partial_\xi^\beta K'(z,x,\xi) \leq
    C (1+|\xi|)^\mu (1+(z-z_0)(1+|\xi|)^\gamma)^{\kappa+
        \frac{|\alpha''|+|\beta''|+|\alphatilde|}{L}} (1+I)^p\\
  \times (1+|\xi|)^\gamma
    \int_{z_0}^z \big((1+|\xi|)^{-\gamma} (1+b)\big)^%
        {1 - \frac{|\alpha'|+|\alphatilde|+|\beta'| +2j}{L}} \, \dd z'.
\end{multline}
The integral on the right hand side can be estimated as in
(\ref{eq:inequality_integral_b}).  Since $l = j + |\alphatilde|$
this yields that $K'$ satisfies the assumptions of
Lemma~\ref{lem:I_bound_symbols2} with $\mu$ and $\kappa$ given by
$\mu' = \mu$ and $\kappa' = \kappa + \frac{2l}{L}$.

Since the $W^{(j)}$ are in $S^{-j(1-\frac{2\gamma}{L})}_%
{1-\frac{\gamma}{L},\frac{\gamma}{L}}$, we can assume that also
the asymptotic sum \linebreak[2]
$\widehat{W} = \sum_{j=0}^\infty W^{(j)}$ is in
$S^0_{1-\frac{\gamma}{L},\frac{\gamma}{L}}((]z_0,Z[\times\R^n) \times\R^n)$.
Then by (\ref{eq:estimate_difference_E_W}) the symbol of $W$ is also in
$S^0_{1-\frac{\gamma}{L},\frac{\gamma}{L}}((]z_0,Z[\times\R^n)\times\R^n)$.
By the Egorov theorem and (\ref{eq:z_der_Btilde}) it follows that
$\widetilde{W}$ is also in $S^0_{1-\frac{\gamma}{L},\frac{\gamma}{L}}%
((]z_0,Z[\times\R^n)\times\R^n)$.

Let $\psi_1 \in C_0^\infty(]z_0,Z[)$ be a function of $z$ only. It is
sufficient to show that for each such $\psi_1$, the operator
$\psi_1 E$ is a Fourier integral operator. Let
$\psi_2 \in C_0^\infty(]z_0,Z[)$ be $1$ on $\supp(\psi_1)$.
Let $\chi = \chi(D_z,D_x)$ be in $\Op S^0(\R^{n+1} \times \R^{n+1})$,
with symbol $\chi(\zeta,\xi)$ that is $1$ on a small conic
neighborhood of $\xi =0, \zeta \neq 0$, and $0$ for $(\zeta,\xi)$,
with $|(\zeta,\xi)| \geq 1$ and outside a larger conic neighborhood
of $\xi =0, \zeta \neq 0$. We assume in particular that $\chi = 0$
on a neighborhood of the set $\zeta = a(z,x,\xi)$. Since $\psi_2$
is only a function of $z$ it commutes with $W$, and we have
\begin{equation}
  \psi_1 E
  = \psi_1 W ((1-\chi) + \chi) \psi_2 E_0 .
\end{equation}
The operator $\chi \psi_2 E_0$ is a smoothing operator.
The proof of Theorem~18.1.35 of \cite{hormander85a} shows that
$\psi_1 W (1-\chi)$ is a pseudodifferential operator in
$\Op S^0_{1-\frac{\gamma}{L},\frac{\gamma}{L}}(\R^{n+1} \times \R^{n+1})$.
Since $E_0$ is a Fourier integral operator in
$I_{1}^{-1/4}(]z_0,Z[ \times \R^n, \R^n; C_0)$, it follows that
$\psi_1 W ((1-\chi) \psi_2 E_0$ is a Fourier integral operator in
$I_{1-\frac{\gamma}{L}}^{-1/4}(]z_0,Z[ \times \R^n, \R^n; C_0)$.
Thus $\psi_1 E$ is a Fourier integral operator.
\end{proof}

\section{Symbols $b$ satisfying the assumption}

We discuss a class of examples such that (\ref{eq:assumption_b1}) holds
for $b$. Let $h$ be a scalar function
$h: \R \rightarrow [0,\infty[ : y \mapsto h(y)$ that satisfies
\begin{equation} \label{eq:h_basic_assumption}
  \text{$h$ is $C^\infty$; $h(y) = 0$, for $y\leq 0$; $h(y) > 0$ for $y>0$} .
\end{equation}
We assume that for $|\xi| > 1$, $b$ is given by
\begin{equation} \label{eq:define_b_general}
  b(z,x,\xi) = |\xi|^\gamma W(z,x,\tfrac{\xi}{|\xi|})
    h(\rho(z,x,\tfrac{\xi}{|\xi|})) ,
\end{equation}
$|\xi| > 1$, where $\rho$ and $W$ are real valued $C^\infty$ functions.

For $h$ we assume that there are an interval $[-\beta,\beta]$ an integer $L$
and a constant $C$ such that on $[-\beta,\beta]$
\begin{equation} \label{eq:h_derivative_alpha_bounds}
  \left| \dede{^j h}{y^j}(y) \right| < C h(y)^{1-j/L},
    j=1,\ldots, L-1  .
\end{equation}
We have in mind the well known example
\begin{equation} \label{eq:h_choice}
  h(y) = \left\{ \begin{array}{ll}
        0          & y \leq 0 , \\
        \exp(-1/y) & y > 0 . \end{array} \right .
\end{equation}
It can be seen that this choice of $h$ satisfies the properties
(\ref{eq:h_basic_assumption}) and (\ref{eq:h_derivative_alpha_bounds}) by
computing the successive derivatives and using that $y^{-j}
\exp(-(1-\alpha)/y)$ is bounded for each $j$, $\alpha<1$ and $y>0$.
By computing the derivatives of $b$ and the assumption
(\ref{eq:h_derivative_alpha_bounds}) we have the following proposition.

\begin{proposition}
Let $h$ be given by (\ref{eq:h_choice}) or otherwise let $h$ satisfy
(\ref{eq:h_basic_assumption}) and suppose that there are $\beta>0$, an
integer $L$ and a constant $C$ such that on $[-\beta,\beta]$
the inequalities (\ref{eq:h_derivative_alpha_bounds}) hold.
Suppose $\rho = \rho(z,x,\frac{\xi}{|\xi|})$
and $W = W(z,x,\frac{\xi}{|\xi|})$ are real valued $C^\infty$ function such
that
\begin{equation} \label{eq:W_general_positive}
  W(z,x,\tfrac{\xi}{|\xi|}) \geq \text{ constant } > 0 .
\end{equation}
Let $b$ satisfy (\ref{eq:define_b_general}) for $|\xi| > 1$.
Then $b$ satisfies (\ref{eq:assumption_b1}).
\end{proposition}

\subsection*{Acknowledgements}

I would like to thank J.~Sj\"ostrand for useful suggestions. The
European Union is acknowledged for its financial support through the
Marie Curie program.

%\bibliographystyle{plain}
%\bibliography{../../Bibfiles/research,../../Bibfiles/research3,%
%  ../../Bibfiles/research2}

\def\cprime{$'$}

\end{document}